\documentclass[14pt]{article}
\usepackage[utf8]{inputenc} 
\usepackage{amsfonts} 
\usepackage{amssymb,amsmath} 
\usepackage{amsthm} 
\usepackage[english]{babel}
\usepackage{euscript} 
\usepackage{bm} 

\usepackage{indentfirst}
\usepackage{pb-diagram}
\newtheorem{lemm}{Lemma}
\newtheorem{theorem}{Theorem}

\newtheorem{remark}{Remark}
\newtheorem{Relation}{Corollary}
\newtheorem{PREDLOZHENIE}{Proposition}

\usepackage[pdftex]{graphicx}
\graphicspath{{Rever/}}
\usepackage{epstopdf}
\usepackage[T2A]{fontenc}

\begin{document}
	
 \date{} 
	
\title{\textbf{LACUNAS AND RAMIFICATION OF VOLUME FUNCTIONS AT SIMPLE ASYMPTOTIC HYPERPLANES AND MONODROMY OF BOUNDARY FUNCTION SINGULARITIES}}
\author{N. ARTEMOV}

\maketitle

\textbf{\abstractname.} 

The {\em volume function} defined by a domain in Euclidean space $\mathbb{R}^n$ is the function on the space of affine hyperplanes equal to volumes cut by these hyperplanes from the domain. The study of these functions originates from the works of Archimedes and Newton (see [8]) and is closely related to the theory of lacunas of hyperbolic partial differential equations.

The volume functions are regular at the hyperplanes of general position with respect to the boundary of the cut domain. We study their behavior at the non-regular hyperplanes, which are either tangent to the boundary of the domain at its finite points or have asymptotic direction. In both cases the local regularity of the restriction of the volume function to a local connected component of the set of regular planes depends on the triviality of a certain relative homology class, the (generalized) even Petrovsky class. In the first case (of finite tangencies) the study of these classes and enumeration of components of regularity (so called local lacunas) at simple singularities of wave fronts is essentially done in [1]; it is formulated in terms of deformations of simple real function singularities. 

We show that the analogous study for asymptotic hyperplanes is related in the same way with the study of boundary function singularities introduced in [7]. We define and calculate local Petrovsky classes for this case and find lacunas, i.e. the local components of the complement of corresponding discriminant in which all local Petrovsky classes are trivial. If component is not lacuna we prove that the corresponding volume function cannot be algebraic. Also we calculate the local monodromy groups, describing the ramification of volume functions.

	\subsection*{Introduction}
Problems of lacunas near singularities of differentiable functions and related problems of the local algebraicity of the volume function appeared as a result of reformulation of problems of lacunas of hyperbolic operators and on the behavior of solutions of hyperbolic PDE's near wave fronts in terms of singularity theory and integral geometry. V. A. Vassiliev in the monograph [1] enumerates local lacunas near simple and many unimodal ordinary real function singularities. In this paper we define and study lacunas near simple boundary real singularities $B_{k}, C_{k}, F_{4}$, where the role of boundary is played by the set where integration form is non-singular (in the case of volume function it is just the <<infinity>> hyperplane in $\mathbb{RP}^{n}$) . Our study is large based on the results of M. A. Gudiev [2], who has enumerated all complements of the corresponding real discriminants. Local monodromy groups are also calculated and results are obtained on the local algebraicity of volumes cut off by hyperplanes from the surfaces, which is graphs of the corresponding morsifications near asymptotic hyperplanes, using obstructions to local algebraicity obtained by V. A. Vasiliev (see [1], [3]).

One of the main problems in the theory of lacunas is the calculation of local Petrovsky classes (in particular, testing them for nontriviality). Consider a miniversal deformation of a smooth function singularity $f$ (see [1], p. 47). The set of deformation parameters ${\lambda}$, such that the corresponding function $f_{\lambda}$ has a critical point on the zero level set $V_{\lambda}$ is called the \textit{discriminant} of the deformation and is denoted by $\Sigma$. Considering real deformations of real germs and real critical points, we obtain the definition of the real discriminant $\Sigma_{\mathbb{R}}.$ Generally speaking, $\Sigma_{\mathbb{R}}$ does not coincide with the real part of $\Sigma,$ since the function $f_{\lambda}$ can have complex critical points on the zero level set. An even local Petrovsky cycle $\text{П}_{ev}(\lambda)$ is an oriented set of real points $\text{Re}V_\lambda$ of a typical local level set $V_\lambda=f_\lambda^{-1}(0)$ real isolated singularity $f:(\mathbb{C}^n,\mathbb{R}^n,0)\rightarrow(\mathbb{C},\mathbb{R},0)$ at zero. A local component of the complement of the real discriminant $\Sigma_{\mathbb{R}}$ of a versal deformation of the germ $f$ is called an even formal lacuna if for all $\lambda$ close to zero from this component the class of cycle $\text{П}_{ev}(\lambda)$ is trivial in the group $H_{n-1}(V_\lambda,\partial V_\lambda)$. This property guarantees the regularity of the analytic continuation of the volume function from the containing $\lambda$ component of the complement to $\Sigma$ to the neighborhood of the origin in the deformation space.

We define lacunas and cycles $\text{П}(\lambda)$ in a similar way for the cases of boundary function singularities $B_{k},C_{k},F_{4}$, in which an improper hyperplane plays the role of a boundary. We will consider sets $V_\lambda\backslash S$ instead of sets $V_{\lambda}$ (where $S$ is the distinguished hyperplane, corresponding to the boundary). All lacunas for the considered boundary singularities and their stabilizations are enumerated in Theorem 1. 
	
	Graphs of morsifications $f_{\lambda}$ (and their stabilizations), which form a locally trivial bundle over the set of typical points of a \textit{generic} line $L\subset(\mathbb{CP}^{n+1})^*$, are smooth hypersurfaces $A_{L}\subset \mathbb{C}^{n+1}$. Let $D\subset\mathbb{C}^{n+1}$ is a ball, which contains the critical points $f_{\lambda}$ and $f_{\lambda}|_{S}, A=D\cap A_{L}.$ For the set $A\cap\mathbb{R}^{n+1}\backslash S$ let us formulate a local variation on the theme of a well-known problem dating back to I. Newton (see [8], [1]):
	
	 can the volume function cut off from $A\cap\mathbb{R}^{n+1}\backslash S$ by the hyperplane $X$  coincide with the algebraic function on connected components $(\mathbb{RP}^{n+1})^*\backslash\Sigma_{\mathbb{R}}$?
	
    The answers to this question are given in Theorem 2. 
	
	In the second part of the paper we calculate for boundary function singulariries $B_{k}, C_{k}, F_{4}$ the homology groups  $H_{n-1}(A\backslash S,A\cap X\backslash S)$ of complex analytic sets $A\backslash S$ reduced modulo typical hyperplane sections $A\cap X$ (see Theorem 3). This groups are generated by vanishing <<caps>> ([1], p. 119) and \textit{Leray tubes} ([1], p. 35). We also calculate the monodromy groups (see Lemmas 2-4), acting on these homology and, in particular, controlling the ramification of the compact parts of the volumes to be cut off.

	 \subsection*{Notations and definitions}
	 Simple boundary function singularities were classified by V. I. Arnol'd [7]. Their real normal forms and miniversal deformations are given below (the cases where $f$ does not have a critical point at zero, corresponding to ordinary simple singularities on the boundary, are not considered):

	 	\begin{center}
	 		
	 		\begin{tabular}{|p{1cm}|p{2cm}|p{8cm}| }\hline
	 			$B_{k}^{\pm}$ & $x^{k}\pm y^2$ &$x^{k}\pm y^2+\lambda_{k-1}x^{k-1}+\lambda_{k-2}x^{k-2}+\dots+\lambda_{1}x+\lambda_{0}$ \\\hline
	 			$C_{k}^{\pm}$ &$xy\pm y^{k}$ &$xy\pm y^{k}+\lambda_{k-1}y^{k-1}+\lambda_{k-2}y^{k-2}+\dots+\lambda_{1}y+\lambda_{0}$ \\\hline
	 			$F_{4}^\pm$ &$\pm x^2+y^3$ &$\pm x^2+y^3+\lambda_{3}x+\lambda_{2}y+\lambda_{1}xy+\lambda_{0}$ \\\hline

	 		\end{tabular}
	 		
	 	\end{center} 
In all cases the boundary $S$ is given by the equation $x=0.$ Along with the singularities of $f$, we will also consider their $(r,s)-$stabilizations $F^{r,s}=f+Q,$ where $Q=\sum_{i=1}^{r}z_{i}^{2}-\sum_{j=r+1}^{r+s}z_{j}^{2}$ is any non-degenerate quadratic form of additional arguments.

	 Let $f_{\lambda}$ be some morsification of the boundary (or ordinary) singularity $f$ (or its stabilization), corresponding to the nondiscriminant point $\lambda$. Consider the line $L=\{f_{\lambda}-c,c\in\mathbb{C}\}$ in the base of the miniversal deformation $f$. For sufficiently small $\lambda$ the line $L$ is generic and intersecting the discriminant $\Sigma$ transversally at $\mu(f)$ points (in the case boundary function singularities $\mu(f)$ is sum of multiplicities of critical points at origin for $f$ and $f|_{S}$, and $\Sigma$ consists of two hypersurfaces in the deformation space: the level set can be either non-smooth or non-transversal to the boundary). Take $c$ as a local coordinate on $L$. Define the hypersurface $A_{L}$ in the space $\mathbb{C}^{n+1}$ with coordinates $(x,y,z_{1},z_ {2},\dots,z_{n-2},c)$ as the graph of $f_{\lambda}$. Its section by the hyperplane $X_{c'}=\{c=c'\}$ coincides with the level set $f_{\lambda}=c'.$ Denote by $D\subset\mathbb{C}^{n+1}$ a disk, containing all critical points of $f_{\lambda}$ together with some neighborhoods, by $A$ the set $A_{L}\cap D$ and by $V_{\lambda}$ the set $A\cap X_{0}$, which is local zero level set of $f_{\lambda}$. Let $F'=f_{\lambda}-c_{1}z_{1}-c_{2}z_{2}-\dots- c_{n-2}z_{n-2}-c_{n-1 }x-c_{n}y-c$ be a linear deformation of $f_{\lambda},\Sigma_{F'}$ be its discriminant. The deformation parameters $(c_{1},\dots,c_{n})$ define local affine coordinates near the point $X_{0}\in(\mathbb{CP}^{n+1})^*,$ where $ (\mathbb{CP}^{n+1})^*$ is the \textit{projective dual} space, consisting points, which are dual to the projective closures of the hyperplanes $X\subset\mathbb{C}^{n+1}.$ By the local variant Zariski's theorem (see [4]) the simple loops around points of $L\cap\Sigma_{F'}$ generate the fundamental group $\pi_{1}(\mathbb{CP}^{n+1})^*\backslash\Sigma_{F '}$.

	 Now we consider only the case of the boundary singularity $f$. Denote by $S\subset\mathbb{C}^{n+1}$ the improper hyperplane given by the equation $\{x=0\}$ (thus, $S\cap X_{c'}$ coincide with the boundary and there will be no contradictions with the notation $S$ introduced earlier). Consider the chain of homomorphisms
	 \begin{equation}\label{1}
	H_{n+1}(D\backslash S,A\cup X\cup\partial D\backslash S)\rightarrow H_{n}(A\backslash S,A\cap X\cup\partial D\backslash S)\rightarrow H_{n-1}(A\cap X\backslash S,\partial D\backslash S).
	 \end{equation}

	The first homomorphism is defined as a composition of mappings from suitable exact sequences of triples and an excision isomorphism, the second is a relative version of the \textit{Mayer-Vietoris differential} (see [1], p. 84). By Thom's isotopy lemma (see [4]), for all typical $X\in(\mathbb{CP}^{n+1})^*\backslash\Sigma_{F'}$ sufficiently close to $X_{0}$, all group types considered in the chain \eqref{1} are isomorphic to each other, being identified along the path segments $\gamma\subset(\mathbb{CP}^{n+1})^*\backslash\Sigma_{F'}$ using \textit{Gauss-Manin connection}. The group $H_{n+1}(D\backslash S,A\cup X\cup\partial D\backslash S)$ contains an distinguished element $\Omega$, corresponding to the oriented union of sets that do not intersect $S$ and are cut off by the hyperplane $X$ from $\text{Re}A$ (if the connected component $\text{Re}A$ which contains such a set does not intersect $S$, then choose any of the two cut off by $X$ parts). Let $\tilde{\text{П}}\in H_{n-1}(A\cap X\backslash S,\partial D\backslash S)$ is the image of the cycle $\Omega$ under homomorphism composition \eqref {1}. Denote by $\tilde{\text{П}}$ the set  $\cup_{i}\tilde{\text{П}}_{i},$ where $\tilde{\text{П}}_{i}$ are connected components of the set $\tilde{\text{П}},\Omega_{i}$ are the corresponding elements of the group $H_{n+1}(D\backslash S,A\cup X\cup\partial D \backslash S)$. Let's define $\textit{volume function}$ $V_{i}$ as follows:
	\begin{equation}\label{2}
	V_{i}=\int_{\Omega_{i}(X)} \omega,
	\end{equation}
	where $\omega$ is the volume form, in our local affine coordinates $\omega=\dfrac{1}{x}dx\wedge dy\wedge dz_{1}\wedge \dots\wedge dz_{n-2}$ . 
	
Now we are ready to formulate the problem: can the volume function $V_{i}$ coincide with some algebraic function on $(\mathbb{RP}^{n+1})^*\backslash\Sigma_{\mathbb{R}(F ')}$ at least in some neighborhood of the point $X_{0}$?

  We will show that in all cases under consideration for $n>2$ $$H_{n+1}(D\backslash S,A\cup X\backslash S)\simeq H_{n}(A\backslash S,A\cap X\backslash S)\simeq H_{n-1}(A\cap X\backslash S)$$ and describe the action of generators $\pi_{1}(L\backslash\Sigma)$ on these homology groups with coefficients in the field $\mathbb{Q},$ obtaining a description of local monodromy groups.

Denote by $\text{П}$ the class of the cycle $\text{Re}V_{\lambda}$ in the group $H_{n-1}(A\cap X,\partial D)$, let $\text{П}=\cup_{i}\text{П}_{i}$ be the decomposition to connected components. The local component $\text{Г}$ of the complement to the real discriminant $\Sigma_{\mathbb{R}}$ of the boundary function singularity $f$ is called \textit{lacuna} if for any $\lambda\in\text{Г}$ any cycle $\text{П}_{i}(\lambda)$  intersects the boundary or is homologous to zero in the group $H_{n-1}(V_{\lambda}\backslash S,\partial V_{\lambda }\backslash S) $.

\subsection*{Lacunas and algebraicity}

\begin{lemm} Each connected component $\text{П}_{i}$ not intersects with boundary of the zero real local level set $\text{Re}V_{\lambda}$ for simple boundary function singularities $f$ of two variables for typical near zero real $\lambda$ realizes non trivial gomology class $\text{П}_{i}(\lambda)\in H_{1}(V_{\lambda}\backslash S,\partial V_{\lambda}\backslash S).$

The same is true for all possible stabilizations $f$ and classes\\ $\text{П}_{i}(\lambda)\in H_{1+r+s}(V_{\lambda}^{r,s}\backslash S,\partial V_{\lambda}^{r,s}\backslash S)$.
\end{lemm}

	\begin{figure}\caption{Examples of typical level sets for $B_{k}$ (left) and $C_{k}$ (right). The boundary $S$ is a vertical line.}
		\center{\includegraphics[width=1\linewidth]{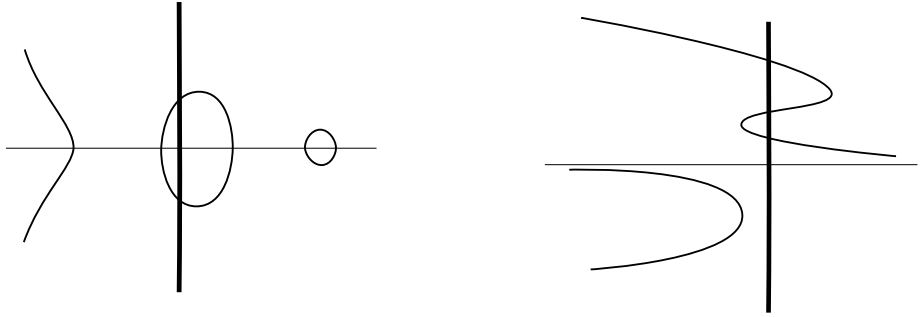}}
		\label{F8}
	\end{figure} 
	
\textbf{Proof.}
In the case of the $B$ series, each component $\text{Г}_{p,q}$ of the complement to the real discriminant is characterized by a pair of numbers $(p,q)$, namely the number of positive and negative points of intersection of the level set $V_\lambda ,\lambda\in\text{Г}_{p,q}$ with line $y=0$ (see [2] and see fig. 1). If $\text{П}_{i}$ is a compact component represented by a flat oval, then $\text{П}_{i}\cap\{y=0\}=\{a,b\} .$ Consider a path in $\text{Г}_{p,q}$ starting at the point $\lambda$ and ending at a typical point of the discriminant, along which the points $a$ and $b$ are approaching each other on the real line and collide at a nonboundary point $p$ at the end point $\lambda'$ of the path, while the other intersection points $\text{П}\cap\{y=0\}$ are fixed. Near $\lambda'$ the discriminant is smooth and its complement consists of two components, one of which is $\text{Г}_{p,q}.$ Then, in the homology $H_{1}(V_\lambda)$ of a typical fiber of the surface $A_{L},$ of the corresponding line $L=\{f_{\lambda'}-c\}$, $\text{П}_{i}$ realizes the basiс vanishing cycle $\Delta$. Due to the connectedness of the Dynkin graph of the singularity $A_{k-1}$ for which any point of the line $L$ close to the discriminant realizes a Morse perturbation, for $k>1$  there exists another basic vanishing cycle $\Delta'$ with $\langle \Delta,\Delta'\rangle\neq0,$ which implies the non-triviality of the cycle $\Delta\in H_{1}(V_\lambda,\partial V_{\lambda})$ and its class $\text{П}_{ i}$ in the group $H_{1}(V_\lambda\backslash S,\partial V_{\lambda}\backslash S)$ (see Proposition 3, p. 12). The case $B_{2}^+$, corresponding to $A_{1}$, is trivial and will be consider separately: a typical fiber is in this case homeomorphic to a cylinder without two points that lie on different sheets of the corresponding two-sheeted Riemann surface and prevent the only vanishing cycle from being pushed out to the boundary.

If the component $\text{П}_{i}$ is not compact, we reason as follows. Let $\text{П}_{i}\cap\{y=0\}=a.$ If there are other real points in the set $\text{П}\cap\{y=0\},$ then consider such point $b$ closest to $a$. Then, colliding the points $a$ and $b$ as above (possibly bypassing the point of the boundary component of the discriminant in the complex domain), we find the vanishing cycle $\Delta$; it easy to see that $\langle\text{П}_{i}, \Delta\rangle=\pm1.$ Hence, due to Poincaré duality $\text{П}_{i}$ realizes non-trivial class in both groups $H_{1}(V_\lambda,\partial V_{\lambda}) $ and $H_{1}(V_\lambda\backslash S,\partial V_{\lambda}\backslash S).$ If $\{a\}$ is the only real point in $\text{П}\cap\{y= 0\}$, then a pair of complex conjugate points always exists in $\text{П}\cap\{y= 0\}$. In this case, it suffices to consider the path in the base, along which this pair of points at first collides at real critical point $b$ with $|b|<|a|$, while the other points of the set are fixed, and then again the points $a$ and $ b$ collide.

 In the case of the $C$ series, every typical fiber $V_{\lambda}$ consists of two non-compact components (see [2] and see fig. 1). On any typical line $L=\{f_{\lambda}-c\}$ there will be a point $\lambda',$ in the fiber over which the cycle $\Delta\in H_{1}(V_\lambda)$ vanishes at the point of transversal intersection of two components. Then $\langle \text{П}_{i},\Delta \rangle=\pm1$, whence $\text{П}_{i}$ realizes non-trivial class in group $H_{1}(V_\lambda\backslash S,\partial V_{\lambda}\backslash S).$
 
 In the case $F_{4}$ the group $H_{1}(V_\lambda)$ is generated by two vanishing cycles with intersection index $\pm1$, one of which  is the oval under consideration (see [2] and see fig. 2).
 
For all possible $(r,s)-$stabilizations $F^{r,s}=f+Q,$ where $r,s$ is inertia indexes of non-degenerate quadratic form $Q,$ the critical values of $f$ and $F$ coincide and there is a natural one-to-one correspondence between the components of the complements to the discriminants $f$ and $F^{r,s}.$ It follows from Gabrielov's formulas [6] that if $|\langle\Delta,\Delta'\rangle|=\kappa$ for vanishing cycles in the homology of the nonsingular fiber $f_{\lambda} $, then for the stabilizations $\Delta^{r,s},\Delta'^{r,s}$ of this cycles in corresponding basis for any stabilization $F_{r,s}$ hold $|\langle\Delta^{ r,s},\Delta'^{r,s}\rangle|=\kappa$ (the same is true for the intersection of absolute and relative cycles, since computation is locally), which implies the second statement of the lemma.  
 $\blacksquare$ 
 
 \begin{Relation}
 	In all cases considered $\text{Г}$ is a lacuna iff for any $\lambda\in\text{Г}$ $\text{Re}V_{\lambda}=\varnothing$ or all connected components $\text{П}_{i}(\lambda)$ of the sets $V_{\lambda}$ intersect the boundary $S$. 
 \end{Relation}

\begin{theorem}The number of lacunas in the space of versal deformation of the considered boundary function singularities of two variables is expressed by the following table:
	\begin{center}
		
		\begin{tabular}{|p{2cm}|p{2cm}|p{2cm}| }\hline
		Title & Number of components & Number of lacunas \\\hline
			$B_{2k}^{+}$ &$(k+1)^2$ &$2$ \\\hline
			$B_{2k}^{-}$ &$(k+1)^2$ &$1$ \\\hline
			$B_{2k+1}^{\pm}$ &$(k+1)(k+2)$ &$1$ \\\hline  
		$C_{2k}^{\pm}$ &$(k+1)^2$ &$k^2$ \\\hline 	$C_{2k+1}^{\pm}$ &$(k+1)(k+2)$ &$k^2+k$ \\\hline $F_{4}^\pm$ &$8$ &$4$ \\\hline

		\end{tabular}
		
	\end{center} 
	
	In the cases $B_{k}^+, F_{4}^+,C_{2k}^+$ (respectively $B_{k}^-, F_{4}^-,C_{2k}^-$) for all stabilizations with $i_{-}=0$ (respectively, $i_{+}=0$) the answer is the same as for $n=2$ ($i_{-},i_{+}$ is respectively the negative and positive indices of inertia of the non-degenerate quadratic form  $Q$). In the cases $C_{2k+1}^{\pm}$ for stabilizations with $i_{+}=0$ or $i_{-}=0$ the number of lacunas equals $(k+1)^2.$  For all other stabilizations all components of the complement to $\Sigma_{\mathbb{R}}$ are lacunas.

\end{theorem}
\textbf{Proof.}
\begin{enumerate}\item \textit{Сase $n=2$.}
	In the paper of M. Gudiev [2], the numbers of components are calculated and the topology of typical $f_{\lambda}$ level sets is described. 
	
	In the case of singularities of the $B$ series, it suffices to consider the trivial case $k=1,$ since with the next increase in $k$ by one, the set of typical level sets is replenished with sets topologically different from the old ones only by new connected components.
	
In the case of $C$, all typical level sets consist of two connected components lying above and below the line $y=0$, respectively. Since each component of the complement to the discriminant is characterized by a pair of numbers $(p,q)$, i.e. the number of positive and negative intersection points of the level set with the boundary, with the next increase in $k$ by one, the set of lacunas will be replenished with new ones in the number of possible $(p,q)$-partitions with $pq\neq0$ of the new Milnor number.

Analysis $F_{4}$ is trivial: see the pictures of all level sets in [2], exactly 4 of them have one of their two components an oval that does not intersect $S$. 
	
	\item\textit{Stabilizations}  
	
In the cases $B_{k}$ with non-empty $\text{Re}V_{\lambda}$, consider the Morse function on the set $\text{Re}V_{\lambda}$, defined by the projection onto the $x$ axis. The critical points of this function are $\text{Re}V_{\lambda}\cap\{y=z_{1}=\dots=z_{n-2}=0\}$. If both indices of inertia of the quadratic form $\pm y^2+Q$ are non-zero, then both Morse indices of each critical point are also non-zero, and the separatrices starting from the critical points set up the linearly connection for the set $\text{Re}V_{\lambda}$, and obviously $\text{Re}V_{\lambda}\cap S\neq\varnothing.$ Thus, all components $\text{Г}$ become lacunas. If $i_{+}$ or $i_{-}$ is equal to zero, then the situation is completely analogous to the case $n=2$. If $\text{Re}V_{\lambda}=\varnothing,$ which corresponds to the equation $y^2=p(x),$ where $p(x)$ is a polynomial without real roots, then for any stabilization all sections $\{x=\text{const}\}$ of the set $\text{Re}V_{\lambda}$ have the same topological type, whence the required statement follows.

For $F_{4}$ we reason in the same way, only we consider the projection onto the $y$ axis.

	Finally, consider the case of the $C$ series. If $\text{Г}$ was a lacuna for $n=2$, then $\text{Г}$ will also be a lacuna for any $(r,s)$-stabilization, since $\text{Re}\;V_{\lambda}^{r,s}$ will have 1 or 2 connected components (to prove this consider the projection onto the $y$ axis), and each such component intersects the boundary $S$, since this holds for $n=2$. Assume that $\text{Г}_{p,q}$ is not lacuna for $n=2$. In the case $C_{2k}^{+}$ $pq=0$, $p+q$ is even and corresponding non-trivial component $\text{П}_{i}\subset\{y<0\}(\text{or}\;\{y>0\})$. Then for stabilization with $i_{-}=0$ the set $\text{Re}\;V_{\lambda}\cap\{y=0\}$ is given by equation $\sum_{i=1}^{n-2}z_{i}^2+\alpha=0$, where $\alpha>0$, sequently the set $\text{Re}\;V_{\lambda}$ have two connected component. Corresponding $\text{П}_{i}$ component again not intersects with boundary $S,$ since the intersection $\text{Re}\;V_{\lambda}\cap\{x=0\}$ is given by equation $\sum_{i=1}^{n-2}z_{i}^2+p(y)=0,$ and $p(y)>0$ for $y<0(\text{or}\;y>0).$ For stabilization with $i_{-}\neq0$ the set $\text{Re}\;V_{\lambda}$ for all typical $\lambda$ is connected and all components $\text{Г}$ realize a lacunas. In the case $C_{2k+1}^{\pm}$ if $\text{Г}_{p,q}$ is not lacuna then $pq=0$, $p+q$ is odd and number of non-lacunas equals $2(k+1).$ For stabilization with $i_{-}=0$ (or $i_{+}=0$), reasoning as above, we conclude that exactly one of two components $\text{Г}_{0,t}$ and $\text{Г}_{t,0}$ for all admissible $t$ becomes lacuna, and the other remains non-lacuna. The analysis of other variants is similar.$\blacksquare$ 

\end{enumerate}

\begin{theorem} If the component of complement of discriminant is not a lacuna, then the volume function \eqref{2} defined above cannot be algebraic.
\end{theorem}
\textbf{Proof.}
	\begin{enumerate}
		\item \textit{Case of even $n$}. The proof is similar to the proof of Theorem 1 in [3]. By virtue of Lemma 1 and Poincaré duality, in the homology group $H_{n-1}(A\cap X)$ there is a basic vanishing cycle $\Delta$ having a nonzero intersection index with the cycle $\text{П}_{i}$, considering as element of group $H_{n-1}(A\cap X, \partial D)$.
		Consider the path $\gamma\subset L,$ connecting the point zero with the point $\lambda',$ in the fiber above which the cycle $\Delta$ vanishes. Then, it follows from Picard-Lefschetz formula that аcting on $\text{П}_{i}(\lambda)$ by consequent rotations around $\lambda'$ along simple loop, corresponding to the path $\gamma$, after the $k$-th iteration cycle $(-1)^{\frac{n(n+1)}{2}}k\langle \text{П}_{i},\Delta \rangle \Delta$ will be added to cycle $\text{П}_{i}$. Thus, each new travel increases the volume function by a nonzero integral over the element $\delta\in H_{n+1}(D,A\cup X)$ corresponding to the vanishing cycle $\Delta$; since this integral is not equal to zero, the volume function cannot be algebraic.
		
		\item \textit{Case of odd $n$}. By corollary 2 (see next page) the monodromy groups cannot be expanded into a direct sum, which means that there will always be a vanishing cycle $\Delta\in H_{n-1}(A\cap X\backslash S)$, class $\tau\in H_{n-1}(A\cap X\backslash S)$, represented as the Leray tube around some vanishing cycle on $A\cap X\cap S,$  and its corresponding monodromy operator $\gamma_{\tau}$, such that $\gamma_{\tau}(\Delta)=\Delta\pm \kappa\tau,\kappa\in\mathbb{Z}\backslash0$. Now, if $\text{П}_{i}$ is a vanishing cycle, which non gomological $\Delta$ (if gomological, then everething is proven, see below), then, due to the transitivity of the action of the monodromy group on vanishing cycles (see [5]) there is such an operator $\gamma'$, such that $\gamma'(\text{П}_{i})=\pm\Delta.$ In this case, we act on $\text{П}_{i}$ first by $\gamma'$ and then, sequentially applying $\gamma_{\tau},$ we observe the logarithmic ramification of the corresponding volume function, since $n+1$ is even (see [1], p. 133-136). If  $\text{П}_{i}$ corresponds to the component intersecting $\partial D,$ then we reason as follows. If, when approaching the boundary critical value, $\text{П}_{i}$ approaches to the boundary (more precisely: the image of the standart reduction homomorphism $\text{red}_{\tilde{x}}(\text{П}_{i})$ modulo the complement of the small  ball $B(\tilde{x})$ around some boundary critical point $\tilde{x}\in S$ realizes the generator of the group $H_{n-1}(B(\tilde{x})\cap A\cap X\backslash S,\partial B(\tilde{x}))$, then we reason as above. If not, by Poincaré duality there exist vanishing cycles that intersect $\text{П}_{i}$. Let $\Delta_{1},\dots,\Delta_{q}$ be all such cycles. If $\Delta$ is one of them, then at first we act on $\text{П}_{i}$ by simple bypass around the critical value for $\Delta$, then sequentially applying on the result by operator $\gamma_{\tau}$ and reason as above. If no, then consider the shortest path $\Delta\rightarrow\Delta_{q+m}\rightarrow\dots\rightarrow\Delta_{q+2}\rightarrow\Delta_{q+1}\rightarrow\Delta_{r}$ in the Dynkin grapf between vertex $\Delta$ and the set $\{\Delta_{j}\}_{j=1}^{q}$, ending in the vertex $\Delta_{r}$ for some $r\in\{1,\dots,q\}.$ We will sequentially act by simple bypasses around the critical values at which the cycles of this chain vanish, starting with $\Delta_{r}$ and ending by $\Delta$, moving along the chain from right to left. As a result of action of the composition of $k+1$ bypasses on $\text{П}_{i}$, the cycle $\Delta_{q+k}$ will be added with a non-zero coefficient, since each cycle in this chain intersects only with neighboring. Action on the final result by going around the critical value for $\Delta$ again and again, we finish the proof.

	\end{enumerate}	
	
	\subsection*{Local monodromy groups}
	
	In the final part of the paper we study homology groups $H_{n+1}(D\backslash S,A\cup X\backslash S)$ with rational coefficients for surfaces $A=A_{L}\cap D$, i.e. graphs of typical morsifications (and their stabilizations) of boundary singularities $B_{k}, C_{k}, F_{4}$ with boundary, corresponding to the <<infinity>> hyperplane. For all three cases with $n>2$ we establish isomorphisms $H_{n+1}(D\backslash S,A\cup X\backslash S)\simeq H_{n}(A\backslash S,A\cap X \backslash S)\simeq H_{n-1}(A\cap X\backslash S)$, we construct basis for these groups and describe the monodromy groups. By analogy with the ordinary singularities, we will call a basis of $H_{n-1}(A\cap X\backslash S)$ \textit{distinguished} (see [1], [5]), iff it is constructed using two distinguished systems of paths for both $f$ and $f|_{S}$ (construction details are described below in the proof of the Theorem 3). The main results are collected in the Lemmas 2-4 and following theorem:
	
	\begin{theorem}
		For free basis $H_{n+1}(D\backslash S,A\cup X\backslash S)$ we can take a set of cycles such that the Mayer-Vietoris differential $\partial:H_{n+1}(D\backslash S ,A\cup X\backslash S)\rightarrow H_{n-1}(A\cap X\backslash S)$ for $n>2$ maps this set isomorphically to the distinguished basis of $H_{n-1}(A\cap X\backslash S)$, formed by union of classes of vanishing basis cycles  from $H_{n-1} (A\cap X)$ and classes of Leray tubes around the vanishing basis cycles of the group $H_{n-2}(A\cap X\cap S).$

	\end{theorem}
	
	The proof of the theorem is given by propositions 1-6. For $n=2$ the monodromy groups are described (Lemmas 2-4): the actions of operators corresponding to simple loops in $L$ around points $L\cap\Sigma$ on the distinguished basis are calculated.
	
	\begin{Relation}
		For any choice of generic line $L$ in space of versal deformation, the monodromy group, acting on the  corresponding $L$ distinguished basis, does not decompose into a direct sum.
	\end{Relation}
	
	\begin{remark}
All elements of the obtained matrices (see Lemmas 2-4) are either the intersection indices of the vanishing cycles, or the coefficients from the calculation of the reduction homomorphism, which also reduces to a (local) calculation of the intersection indices. Therefore, results for any $n>2$ can be obtained by applying Gabrielov's formulas [6] to the same results for $n=2$.
	\end{remark}
	
We will analyze the in most detailed way case $F_{4}$ and will ommit the similar steps in the cases $B_{k}$ и $C_{k}$. 

	\subsection*{$F_{4}$}

	Let us choose a morsification of the germ $f$ in the form $F=x^2+y^3+\sum_{i=1}^{n-2}z_{i}^{2}-y+x$. This morsification corresponds to the hypersurface $A\subset\mathbb{C}^{n+1}$, given by the equation $w=x^2+y^3+\sum_{i=1}^{n-2}z_{i}^{2}-y+x$, and the family of hyperplanes $X_{c}:\{w=-c,c\in\mathbb{C}\}.$ $X_{\nu_{1}}$ and $X_{\nu_{2}}$ are tangent to $A$ at critical points of $F,$  $X_{\nu_{3}}\cap S$ and $X_{\nu_{4}}\cap S$ are tangent to $A\cap S$ at critical points of $F|_{S}$, all critical points $\{x_{i}\}_{i=1}^{4}$ and corresponding critical values $\{\nu_{i}\}_{i=1}^{4}$ are real and $\nu_{4}<\nu_{2}<0<\nu_{3}<\nu_{1}$ (fig. 2 for case $n=2$).  Let us choose in $\mathbb{C}^{n+1}$ with coordinates $(x,y,z_{1},\dots,z_{n-2},w)$ a disk $D$ centered at origin, containing all critical points of $F$ and $F|_{S}$ together with their small neighborhoods. Fix $X\equiv X_{0}.$
	
	\begin{figure}\caption{Rearrangements the level sets for   $F_{4}$}
		\center{\includegraphics[width=1\linewidth]{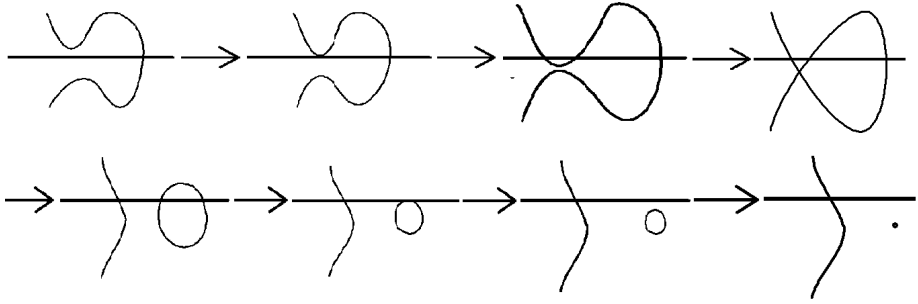}}
		\label{F7}
	\end{figure}
	
	\begin{PREDLOZHENIE}For $n>2$ $H_{n-1}(A\cap X\backslash S)\simeq\mathbb{Q}^{4}$, for $n=2$ $H_{n-1}(A\cap X\backslash S)\simeq\mathbb{Q}^{5}.$
	\end{PREDLOZHENIE}
	\textbf{Proof.}
	Consider a fragment of the exact sequence of the pair $(A\cap X,A\cap X\backslash S):$
	$$0\rightarrow H_{n}(A\cap X,A\cap X\backslash S)\xrightarrow{\varphi_{1}} H_{n-1}(A\cap X\backslash S)\xrightarrow{\varphi_{2}} H_{n-1}(A\cap X)\rightarrow0.$$
    Instead of the first nontrivial term, one can write $H_{n-2}(A\cap X\cap S)$, since
	$H_{n}(A\cap X,A\cap X\backslash S)\simeq H_{n}(N_{\varepsilon},N_{\varepsilon}\backslash S)\simeq H_{n-2}(A\cap X\cap S).$ The first isomorphism here is an excision isomorphism, $N_{\varepsilon}$ denotes the space of the bundle of small disks of the normal bundle to $(A\cap X\cap S)\subset (A\cap X).$ The second isomorphism is Thom isomorphism. Then the homomorphism $\varphi_{1}:H_{n-2}(A\cap X\cap S)\rightarrow H_{n-1}(A\cap X\backslash S)$ is given by the Leray tube operators. In the case $n=2$ the group $H_{0}(A\cap X\cap S)$ is generated by three cycles $a_{0},a_{1},a_{2}$, and the cycles $a_{1}-a_{0},a_ {2}-a_{1}$ vanishing at critical points $x_{3},x_{4}$ of function $F|_{S}$. Denote the images of the cycles $a_{0},a_{1},a_{2}$ under the tube homomorphism $\varphi_{1}$ respetively by $\tau_{0},\tau_{1},\tau_{2}$. As generators of $H_{1}(A\cap X)$ we take the elements $\Delta_{1},\Delta_{2}$ represented by cycles vanishing at the critical points $x_{1},x_{2},$ transported with the help of the Gauss-Manin connection to our fiber along the paths going along the line $\text{Im}=0$ and bypassing the critical values in the half-plane $\text{Im}>0.$ Thus, for the basic generators of $H_{1}(A\cap X\backslash S)$, we can take the elements represented by the cycles $\varphi_{2}^{-1}\Delta_{1},\varphi_{2}^{-1}\Delta_{2},\tau_{0},\tau_{1},\tau_{2}.$ In the case of $n>2$, the group $H_{n-2}(A\cap X\cap S)$ is generated by two vanishing cycles, and a similar argument gives 4 basic generators of the group $H_{n-1}(A\cap X\backslash S).\blacksquare$ 
	
	\begin{PREDLOZHENIE}$H_{n+1}(D\backslash S,A\cup X\backslash S)\simeq \mathbb{Q}^4.$
	\end{PREDLOZHENIE}
	\textbf{Proof.}
	The exact sequence of the triple $(D\backslash S,A\cup X\backslash S,D\cap(X\backslash S))$ implies an isomorphism $H_{n+1}(D\backslash S, A\cup X\backslash S)\simeq H_{n}(A\cup X\backslash S,X\backslash S)$, since $H_{n+1}(D\backslash S,X\backslash S)=H_{n}(D\backslash S,X\backslash S)=0,$ because $X\backslash S$ is a deformation retract of the set $D\backslash S.$ Next, we note that, by the excision property $H_{n}(A\cup X\backslash S,X)\simeq H_{n}(A\backslash S,A\cap X).$ Consider a fragment of the exact sequence of the pair $(A\backslash S, A\cap X\backslash S)$:	
	$$H_{n}(A\backslash S)\rightarrow H_{n}(A\backslash S,A\cap X\backslash S)\xrightarrow{\psi_{1}} H_{n-1}(A\cap X\backslash S)\xrightarrow{\psi_{2}} H_{n-1}(A\backslash S).$$
	$H_{n}(A\backslash S)\simeq H_{n}(\mathbb{C}^{n}\backslash\mathbb{C}^{n-1})=0$, for $n>2$ $H_{n-1}(A\backslash S)\simeq 0$, for $n=2$ $H_{n-1}(A\backslash S)\simeq\mathbb{Q}$ and as a generator one can take a tube around any point of the curve $A\cap S.$ Then for $n=2$ for a basis in $H_{2}(A\backslash S, A\cap X\backslash S)$ we can take elements represented by relative cycles $\{\beta_{i}\}_{i=1}^4$ with boundary $\varphi_{2}^{-1}\Delta_{1},\varphi_{2}^{-1}\Delta_{2},\tau_{1}-\tau_{0},\tau_{2}-\tau_{1},$ since $\psi_{1}$ isomorphically maps $H_{2}(A\backslash S,A\cap X\backslash S)$ onto $\text{Ker}\psi_{2}.$ $\beta_{1}$ and $\beta_{2}$ are realized by <<caps>> ([1], p. 119), $\beta_{3}$ and $\beta_{4}$ are tubes around vanishing segments, defined as elements $H_{1}(A\cap S,A\cap S\cap X)$, connecting respectively $a_{0}$  with $a_{1}$ and $a_{1}$ with $a_{2}$ on the curve $\text{Re}\{A\cap S\}.$ For $n>2$ $\psi_{1}$ is an isomorphism. $\blacksquare$
	
	Fix $n=2$. Let  $\{\gamma_{i}\}_{i=1}^4$ be loops in $D_{0}$ around $\{\nu_{i}\}_{i=1}^4,$ oriented counterclockwise. We orient the cycles $\varphi_{2}^{-1}\Delta_{1},\varphi_{2}^{-1}\Delta_{2}$ so that $\langle\Delta_{1},\Delta_ {2} \rangle=+1$ with respect to the complex orientation of the curve $A\cap X$. The cycles $\tau_{1}-\tau_{0},\tau_{2}-\tau_{1}$ are oriented according to the definition of the tubular operator. These orientations induce orientations of cycles  $\{\beta_{i}\}_{i=1}^4.$
	
	\begin{lemm} The action of the generators $\{\gamma_{i}\}_{i=1}^4$ of the group $\pi_{1}(D_{0}- \cup_{i=1}^{4}{\nu_{i}})$ on the basis $\{\beta_{i}\}_{i=1}^4$ as follows: $\gamma_{i}(\beta_{j})=\beta_{j}+\eta_{ij}\beta_{i},$ where $(\eta_{ij})$ is the following matrix: 
		$$	
		\begin{pmatrix}
		0& 1& 0& 0\\
		-1& 0& 0& 0\\
		1& 0& -2& 1\\
		-1& 1& 1& -2
		\end{pmatrix}.
		$$

	\end{lemm}
	\textbf{Proof.}
	
	The action of $\gamma_{1}$ and $\gamma_{2}$ on $\beta_{1}$ and $\beta_{2}$ is described by generalized Picard-Lefschetz formulas, and on $\beta_{3}$ and $ \beta_{4}$ these operators will have no action: the corresponding boundaries do not intersect.
	
	Under the action of $\gamma_{3}$ on $\beta_{3}$, the following happens: the movement around the loop counterclockwise corresponds to the permutation of the points $a_{1}$ and $a_{2}$ together with counterclockwise rotation of the segment connecting them, and the claim follows from the fact that the monodromy commutes with tubular operation. The action of $\gamma_{4}$ on $\beta_{4}$ is similar. The same reasoning gives answers for the action of $\gamma_{3}$ on $\beta_{4}$ and $\gamma_{4}$ on $\beta_{3}.$
	
	The action of $\gamma_{3}$ on $\beta_{1}$ and $\gamma_{4}$ on $\beta_{2}$ are calculated as follows. Consider a small ball $B_{x_{3}}\subset\mathbb{C}^3$ around the point $x_{3}.$ Since the point $x_{3}$ is not a parabolic point of the divisor $ A\cap S,$ local monodromy is described by Lemmas 7.1 and 7.2 on p. 133 in [1]. The homomorphism $\text{red}$, i.e. the reduction modulo the complement of $B_{x_{3}}$ sends $\Delta_{1}$ to the generator of $H_{1}(B_{x_{3}}\cap A\cap X_{\nu_{3}+\delta}\backslash S,\partial B_{x_{3}}\backslash S),$ and hence $\beta_{1}$ to the generator $H_{2}(B_{x_{3}}\cap A\backslash S,A\cap X\backslash S\cup \partial B_{x_{3}}).$ To calculate the variation over the loop around $\nu_{3}$, we project $X_{\nu_{3}+\delta}\cap B_{x_{3}}\backslash S$ to $S$ and compare the orientations of $\Delta_{1}\cap B_{x_{3}}$ and $\tau_{2}-\tau_{1}.$
	
  The operator $\gamma_{3}$ will not act on $\beta_{2}$, since $\langle \Delta_{1},\Delta_{2}\rangle=1,$ and the result of transportation of $\Delta_{2}$ along the segment $[0,\nu_{3}-\delta]$ does not intersect $B_{x_{3}}:$ it is not linked to the vanishing arc of the oval located above $S$.
	
	It remains to calculate the action of $\gamma_{4}$ on $\beta_{1}$. To do this, we describe the result of transportation of $\Delta_{1}$ along the path from $\nu_{1}-\delta$ to $\nu_{4}+\delta$, bypassing the values $\nu_{3},\nu_{2}$ in the half-plane $\text{Im}>0,$ considering the projection onto $S$ in $X_{\nu}.$ The bypass along the arc in $\text{Im}>0$ of the value $\nu_{3}$ is locally equivalent to a rearrangement of the curve $y=-x^2-c$, where $c$ runs counterclockwise from $\delta$ to $-\delta$, since counterclockwise bypass around $\nu_{3}$ corresponding to the counterclockwise movement  of the points $a_{1}$ and $a_{2}$ in the projection. In the projection onto $S\cap X_{\nu}$, the cycle $\Delta_{1}$ transported in this way will bypass the point $a_{1}$ from below. The bypass around of the critical value $\nu_{2}$ is locally equivalent corresponds to rearrangement in the neighborhood of the critical point $x_{2},$ written as follows: $x^2-y^2=c,$ when $c$ runs counterclockwise from $\delta$ to $-\delta$. Since $\langle\Delta_{1},\Delta_{2}\rangle=1,\text{Var}(\Delta_{1}\cap B_{x_{2}})=-\Delta_{2}.$ $\Delta_{2}$ in $\mathbb{R}^{2}$ is given by $x^2+\tilde{y}^2=\delta,$ where $\tilde{y}=iy.$ Then the ordered frame $(\partial/\partial y,i\partial/\partial y)$ at the point $(\sqrt{\delta},0)$ defines the orientations of the cycles $(\Delta_{1},\Delta_{2}).$ After passing along the arc $\nu_{2}$, considered above, the cycle $\Delta_{2}$ will be written as $\tilde{x}^2-y^2=-\delta,$ where $\tilde{x}=ix.$ At the point $(0,\sqrt{\delta})$ the orienting frame is the pair $\partial/\partial x,i\partial/\partial x$. Since $\text{Var}(\Delta_{1}\cap B_{x_{2}})=-\Delta_{2},$ the part of the cycle $\Delta_{1}$ transported along the arc will go from the point $(0,\sqrt{\delta})$ towards $(0 ,-\sqrt{\delta})$ along an arc in the direction of the positive orientation of $\Delta_{2},$ but $\Delta_{1}$ itself will be oriented in the opposite direction. Generally speaking, the morsification of $F$ can be chosen so that $\nu_{4}$ and $\nu_{2}$ are arbitrarily close (replacing the monomial $z_{1}$ with $\varepsilon z_{1}$), so we can consider the considered fiber close to $X_{\nu_{4}}.$ In the projection, the transported cycle $\Delta_{1}$ will still bypass from below the point $a_{2}$ and go up along the imaginary axis $\text{Im},$ its homology class $H_{1}(B_{x_{4}}\cap A\cap X_{\nu_{4}+\delta}\backslash S,\partial B_{x_{4}}\backslash S)$ equals minus the projection class $\Delta_{2}$. Therefore, the variation of $\beta_{1}$ over $\gamma_{4}$ is equal to $-\beta_{4}$. Note that if we were to bypass $\nu_{3}$ clockwise, when we considered bypass around $\nu_{3}$, the class $\text{red}_{x_{4}}(\Delta_{1})$ would be equal zero, and the $\gamma_{4}$ operator would have no action.
	
	We can check the last calculation as follows: consider two paths $l_{1},l_{2},$ starting from $\nu_{1}$ and going to $0,$ while bypassing $\nu_{3}$ respectively counterclockwise and clockwise. It is clear that $l_{1*}(\Delta_{1})-l_{2*}(\Delta_{1})=-\partial \beta_{3}.$ Now let's act $\gamma_{4}.$ On the one hand, the calculations above imply $\gamma_{4}(l_{1*}(\Delta_{1})-l_{2*}(\Delta_{1}))=\gamma_{4}(l_{1 *}(\Delta_{1}))-\gamma_{4}(l_{2*}(\Delta_{1}))=l_{1*}(\Delta_{1})-\partial\beta_{4 }-l_{2*}(\Delta_{1})=-\partial\beta_{4}-\partial\beta_{3}.$ On the other hand, $\gamma_{4}(l_{1*} (\Delta_{1})-l_{2*}(\Delta_{1}))=\gamma_{4}(-\partial\beta_{3})=-\gamma_{4}(\partial\beta_{ 3})=-\partial\beta_{3}-\partial\beta_{4}.$
	$\blacksquare$
	
	\subsection*{$B_{k}$} 
	
Consider morsification of the corresponding germ $f$, having the form $F=y^2+(x-a_{1})(x-a_{2})\dots(x-a_{k})+\sum_{j=1}^{n-2 }z_{j}^2,$ where all $a_{i}$ are different small positive numbers and $a_{1}<a_{2}<\dots<a_{k}$. In this case function $F$ has $k-1$ critical points $\{x_{i}\}_{i=1}^{k-1}$, and in corresponding singular fibers hold $pr_{x}(x_{i})\in(a_{i},a_{i+1})$. The function $F|_S$ has a single critical point at the origin. Let us denote by $\{\nu_{i}\}_{i=0}^{k}$ the corresponding different (which easy to achieve by a small perturbation of $a_{i}$) critical values. Consider the line in the deformation base $L=\{F-c,c\in\mathbb{C}\}$ and the corresponding hypersurface $A=A_{L}\cap D\subset\mathbb{C}^{n+1}$ in the ball $D.$ Fix a hyperplane $X=\{c=0\}\subset\mathbb{C}^{n+1}.$
	
	\begin{PREDLOZHENIE}$H_{n-1}(A\cap X\backslash S)\simeq\mathbb{Q}^{k}$ for $n>2$, $H_{n-1}(A\cap X\backslash S)\simeq\mathbb{Q}^{k+1}$ for $n=2$.
	\end{PREDLOZHENIE}
	
For distinguished basis of this group we take the vanishing cycles $\{\Delta_{i}\}_{i=1}^{k-1}$ of the group $H_{n-1}(A\cap X)$ corresponding to $A_{k -1}$ and a tube $\Delta_{0}$ around a vanishing cycle in $H_{n-2}(A\cap X\cap S)$ of the boundary section, which is given by the standard Morse equation. The case $n=2$ is not typical: $A\cap X\cap S=\{2\;points\}$.
$\blacksquare$
	
	\begin{PREDLOZHENIE}$H_{n}(D\backslash S,A\cap X\backslash S)\simeq\mathbb{Q}^{k}\blacksquare$.
	\end{PREDLOZHENIE}

	We fix $n=2$.
	Connect the critical values $\{\nu_i\}_{i=0}^{k-1}$ to zero by paths $\gamma_i\subset D_{0}\subset L$ along the $\text{Re}$ axis and bypassing other critical values in the half-plane $\text{Im}>0.$ We orient the vanishing cycles $\Delta_{i}$ (and the corresponding <<caps>> $\beta_{i}$) as follows: if $\Delta_{i}$ is standard oriented oval in $\mathbb{R}^2,$ then cycles $\Delta_{i\pm1}$, that are ovals in $\mathbb{R}^{2}=(\text{Re}\;x,\text{Im}\;y)$, are oriented so that $\langle\Delta_{i},\Delta_{i\pm1}\rangle=+1$ with respect to the complex orientation of the fiber; we orient $\beta_{0}$ so that $\langle\text{pr}_{y}\Delta_{1},\text{pr}_{y}I\rangle=-1$ hold in typical fiber close to $X_ {\nu_{0}}$, where $I$ is a vanishing segment at the boundary around which the tube is taken. We also denote by $\gamma_{i}$ the simple loops around $\nu_{i}$ and the corresponding monodromy operators.
	\begin{lemm}
		Action of generators $\{\gamma_{i}\}_{i=0}^{k-1}$ of the group $\pi_{1}(D_{0}-\cup_{i=0}^{k-1 }{\nu_{i}})$ on the basis elements $\{\beta_{j}\}_{j=0}^{k-1}$ as folows: $\gamma_{i}(\beta_{j}) =\beta_{j}+\eta_{i+1,j+1}\beta_{i},$ where $(\eta_{i+1,j+1})$ is the following matrix:
		$$	
		\begin{pmatrix}
	    -2& 1& 0& .&. &.&0 \\
		0& 0&(-1)^{k}&0 &.& .& 0\\
		.& (-1)^{k-1}& 0&(-1)^{k-1}& 0&.&0 \\
		.& 0& (-1)^{k-2}& .& .&.&0  \\
			.& .&0& .&0& 1&0  \\
		.& .& .&.& -1& 0&-1\\
		0& 0& 0& 0&0& 1& 0
		\end{pmatrix}.
		$$
	\end{lemm}
	     
	 $\blacksquare$
	 
	 	\subsection*{$C_{k}$} 
	 	
	 Consider morsification of the corresponding germ $f$, choosen in the form $F=xy+(y-a_{1})(y-a_{2})\dots(y-a_{k})+\sum_{j=1}^{n-2}z_{ j}^2,$ where all $a_{i}$ are different small positive numbers and $a_{1}<a_{2}<\dots<a_{k}$. We will also additionally assume that $|a_{1}|<<|a_{2}-a_{3}|.$ The function $F|_S$ has $k-1$ critical point $\{x_{i}\} _{i=1}^{k-1}$, and in corresponding singular fibers hold $pr_{y}(x_{i})\in(a_{i},a_{i+1})$. The function $F$ has a single critical point $x_{0}$. Denote by $\{\nu_{i}\}_{i=0}^{k}$ the corresponding critical values. Consider the line in the deformation base $L=\{F-c,c\in\mathbb{C}\}$ and the corresponding hypersurface $A=A_{L}\cap D\in\mathbb{C}^{n+1}$ in the ball $D.$ We fix a hyperplane $X=\{c=0\}\subset\mathbb{C}^{n+1}.$
	 	
	 	\begin{PREDLOZHENIE}$H_{n-1}(A\cap X\backslash S)\simeq\mathbb{Q}^{k}$ for $n>2$, $H_{n-1}(A\cap X\backslash S)\simeq\mathbb{Q}^{k+1}$ for $n=2$.
	 	\end{PREDLOZHENIE}

	For distinguished basis of this group we take the tubes $\{\tau_{i}\}_{i=1}^{k-1}\in H_{n-1}(A\cap X)$ around the vanishing cycles of the group $ H_{n-2}(A\cap X\cap S)$ corresponding to $A_{k-1}$ and the vanishing cycle $\Delta_{0}\in H_{n-1}(A\cap X).$ Case $ n=2$ is atypical: $A\cap X\cap S=\{k\;points\}.\blacksquare$

	 	\begin{PREDLOZHENIE}$H_{n}(D\backslash S,A\cap X\backslash S)\simeq\mathbb{Q}^{k}.\blacksquare$
	 	\end{PREDLOZHENIE}

	 Fix $n=2$.
	 Connect the critical values $\{\nu_i\}_{i=0}^{k-1}$ to zero by paths $\gamma_i\subset D_{0}\subset L$ along the $\text{Re}$ axis and bypassing other critical values in the half-plane $\text{Im}>0.$ Orient the tubes according to the orientation $\{y_{i+1}\}-\{y_{i}\}$ of the vanishing cycles on the boundary. We also denote by $\gamma_{i}$ the simple loops around $\nu_{i}$ and the corresponding monodromy operators.
	 	\begin{lemm}
	 		Action of generators $\{\gamma_{i}\}_{i=0}^{k-1}$ of the group $\pi_{1}(D_{0}- \cup_{i=0}^{k-1 }{\nu_{i}})$ on the basis $\{\beta_{j}\}_{j=0}^{k-1}$ as follows: $\gamma_{i}(\beta_{j}) =\beta_{j}+\eta_{i+1,j+1}\beta_{i},$ where $(\eta_{i+1,j+1})$ is the following matrix:
	 		$$	
	 		\begin{pmatrix}
	 		0& 0& 0& .&. &.&0 \\
	 		1& -2&1&0 &.& .& 0\\
	 		.& 1& -2&1& 0&.&0 \\
	 		.& 0& 1& .& .&.&0  \\
	 		.& .&0& .&-2& 1&0  \\
	 		.& .& .&.& 1& -2&1\\
	 		0& 0& 0& 0&0& 1& -2
	 		\end{pmatrix}.
	 		$$
	 	\end{lemm}
	 
	 For $i,j>0$ the corresponding elements of the matrix are obtained from the fact that the monodromy commutes with the tubular operation. If the critical value $\nu_{0}$ greater than $0$ (for $\nu_{0}<0$ the arguments is similar), then for all others holds $\nu_{2s}>0,\nu_{2s+1}<0$. It is clear that $\text{red}_{x_{i}}(\Delta_{0})\neq0$ exactly when $\Delta_{0}$ is linked to the arc $\varphi_{i}$ between $ a_{i}$ and $a_{i+1}$ in a fiber close to $X_{\nu_{i}}$. It follows from the condition $|a_{1}|<<|a_{2}-a_{3}|$, that in a fiber close to $X_{\nu_{0}}$ (and also close to $X_{\nu_{1}}$ after transporting $\Delta_{0}$ along real line from $\nu_{0}$ to $\nu_{1}$) holds $\langle\Delta_{0},\varphi_{1}\rangle=\pm1.$ It immediately follows that $\text{red}_{x_{1}}(\Delta_{0})\neq0,\text{red}_{x_{i}}(\Delta_{0})=0$ for all $i>1.$ $\blacksquare$

\end{document}